\documentclass[11pt]{article}
\usepackage{graphicx}
\usepackage{pdfsync}
\usepackage{psfrag}
\usepackage{amsmath}
\usepackage{amssymb}
\usepackage{epstopdf}
\usepackage{color}

\definecolor{Red}{rgb}{1.0,0.0,0.0}

\newtheorem{theorem}{Theorem}
\newtheorem{proposition}[theorem]{Proposition}
\newtheorem{lemma}[theorem]{Lemma}
\newtheorem{corollary}[theorem]{Corollary}
\newtheorem{definition}{Definition}

\title{Upper bounds in the Ohtsuki-Riley-Sakuma partial order on 2-bridge knots}
\author{Scott M. Garrabrant, Jim Hoste and Patrick D. Shanahan}

\begin{document}
\maketitle
\begin{abstract}
In this paper we use continued fractions to study a partial order on the set of 2-bridge knots derived from the work of  Ohtsuki, Riley, and Sakuma. We establish necessary and sufficient conditions for any set of 2-bridge knots to have an upper bound with respect to the partial order. Moreover, given any 2-bridge knot $K_1$ we characterize all other 2-bridge knots  $K_2$ such that $\{K_1,K_2\}$ has an upper bound.  As an application we answer a question of Suzuki, showing that there is no upper bound for the set consisting of the trefoil and figure-eight knots.
\end{abstract}

\section{Introduction}
\label{introduction}

Given two knots $K$ and $J$ in $S^3$, an interesting question in knot theory, and one which has received a great deal of attention, is whether there exists an epimorphism between the fundamental groups of the complement of $J$ and the complement of $K$. The existence of such an epimorphism defines a partial order on the set of prime knots. Since the granny knot and the square knot are different composite knots with isomorphic fundamental groups, we see that it is necessary to consider only prime knots for this partial order. Moreover, since there exist different links with homeomorphic complements, this does not define a partial order for links. A related partial order on the set of all knots can be defined by requiring the epimorphism to have the additional property that it preserve peripheral structure, that is, that it take the subgroup generated by the meridian and longitude of $J$ into a conjugate of the subgroup generated by the meridian and longitude of $K$. The reader is referred to work of Silver and Whitten \cite{SW:2006} for more details on this partial ordering.

In this article we focus on 2-bridge knots and a partial ordering $\ge$ associated to a construction given by Ohtsuki, Riley, and Sakuma in \cite{ORS}.  This ordering can be defined in terms of the continued fraction expansions of 2-bridge knots, and this is the point of view that we exploit in this article.  On the other hand, this ordering can also be defined in terms of the existence of a particular type of function, called a branched fold map, between the complements of two 2-bridge knots.  If $J \ge K$, then the branched fold map between the knot complements induces an epimorphism of the fundamental groups which preserves peripheral structure.  Therefore, the Silver-Whitten partial ordering is a refinement of the Ohstuki-Riley-Sakuma partial ordering for 2-bridge knots.  An open question, proposed in \cite{ORS}, is whether these two partial orderings on 2-bridge knots are in fact the same.  Gonzalez and Ramirez \cite{GR:2001} have shown that indeed this is the case if the smaller knot is both a 2-bridge knot and a torus knot.  Additional evidence also appears in a recent paper by Lee and Sakuma \cite{LS:2010}, where it is shown that an epimorphism between 2-bridge knot groups takes upper meridional pairs to upper meridional pairs if and only if it is given by a branched fold map as described in \cite{ORS}.  An affirmative answer to this question would also imply that 2-bridge knots with exactly three distinct boundary slopes are minimal with respect to the Silver-Whitten partial order, and this has recently been proven by the second and third authors \cite{HS:2010}. (Here, a knot $K$ is minimal if whenever $K$ is greater than $J$ in the Silver-Whitten partial order, then $J$ is either $K$ itself or the unknot.)

This paper was motivated by  a question of M. Suzuki: Does there exist a 2-bridge knot $K$ whose group surjects onto the groups of both the trefoil, $K_{1/3}$, and the figure eight knot, $K_{3/5}$? It is known that knots that are not 2-bridge exist with this property. (See the related papers \cite{KT:2008a}, \cite{KT:2008b} and \cite{HKMS:2010}.) In this paper we give a partial answer to Suzuki's question, showing that there does not exist an upper bound of the set $\{ K_{1/3},K_{3/5}\}$ with respect to the Ohtsuki-Riley-Sakuma partial order. Our answer follows trivially from a more general theorem in which we provide a complete classification of all pairs of 2-bridge knots $K_1$ and $K_2$ such that  $\{K_1, K_2\}$ has an upper bound.

The paper will proceed as follows. We begin by recalling basic facts about the classification of 2-bridge knots and develop the necessary tools regarding continued fraction expansions. Next we describe the Ohtsuki-Riley-Sakuma partial ordering in terms of continued fractions. In the final section, we state and prove the main theorem which classifies all pairs of knots having an upper bound in the Ohtsuki-Riley-Sakuma partial order. We close by extending this result to a set of any number of 2-bridge knots.
 
\section{Continued fractions and 2-bridge knots}

Recall that each 2-bridge knot  or link corresponds to a relatively prime pair of integers $p$ and $q$ with $q$ odd in the case of a knot and $q$ even in the case of a  link. We denote the knot or link as $K_{p/q}$. Furthermore, $K_{p/q}$ and $K_{p'/q'}$ are ambient isotopic as unoriented knots or links if and only if $q'=q$ and $p' \equiv p^{\pm 1}\  (\mbox{mod } q)$ (see \cite{BZ:2003} for details).  In this paper, we will not distinguish between a knot or link $K_{p/q}$ and its mirror image $K_{-p/q}$.  Therefore, two 2-bridge knots or links $K_{p/q}$ and $K_{p'/q'}$ are {\it equivalent} if and only if $q'=q$ and either $p' \equiv p^{\pm 1}\  (\mbox{mod } q)$ or $p' \equiv -p^{\pm 1}\  (\mbox{mod } q)$.  We denote the set of equivalence classes of 2-bridge knots and links as $\cal B$.

The theory of 2-bridge knots and links is closely tied to continued fractions. We adopt the convention used in \cite{ORS} and define $r+[a_1, a_2, \dots, a_n]$ as the continued fraction
 
$$r+[a_1, a_2, \dots, a_n]=r+\frac{1}{a_1+\displaystyle \frac{1}
{
\begin{array}{ccc}
a_2+&&\\
&\ddots&\\
&&+\displaystyle \frac{1}{a_n}\\
 \end{array}
 }
 }$$

A given fraction $p/q$ can be expressed as a continued fraction in more than one way. However there are various schemes for producing a canonical expansion. (See for example \cite{BZ:2003}.) These are related to variations of the Euclidean algorithm used for finding the greatest common divisor of two integers.

Given $p$ and $q$, with $q\ne 0$, let $r_0=p$ and $r_1=q$ and then write $r_0$ as $r_0=a_1 r_1+r_2$ where $|r_2|<|r_1|$. If $r_2\ne0$, then repeat the process using $r_1$ and $r_2$, that is, write $r_1$ as $r_1=a_2r_2+r_3$, with $|r_3|<|r_2|$. At each step (except the last) there are exactly two choices for $a_i$: either  $a_i=\lfloor \frac{r_{i-1}}{r_i}\rfloor$ or $a_i=\lceil \frac{r_{i-1}}{r_i}\rceil$, the floor or ceiling, respectively, of $\frac{r_{i-1}}{r_i}$. 
Since the {\it remainders} $r_2, r_3, \dots$ are getting strictly smaller in magnitude, the process must end. This process, which is recorded in the following set of equations,  is known as the {\it Euclidean algorithm} (or, perhaps more properly, a {\it generalized} Euclidean algorithm).
\begin{align*}
\label{euclidean algorithm equations}
r_0&=a_1r_1+r_2\\
r_1&=a_2r_2+r_3\\
&\vdots\\
r_{n-2}&=a_{n-1}r_{n-1}+r_n\\
r_{n-1}&=a_nr_n
\end{align*}

When we arrive at the last pair, $r_{n-1}$ and $r_n$, we have that $r_{n-1}$ is a multiple of $r_n$ and so there is a unique choice for $a_n$. Notice that none of the {\it partial quotients} $a_i$ are zero and also that the last partial quotient, $a_n$, is not $\pm 1$, since $|r_n|<|r_{n-1}|$. It is well known, and easy to show, that the greatest common divisor of $p=r_0$ and $q=r_1$ is $a_n$ and also that $\frac{p}{q}=a_1+[a_2, a_3, \dots, a_n]$.

At each step, one of the two choices for $a_i$ is even and the other odd. We can make the construction deterministic by always making the even choice. Since we have no choice at the last step, this may not produce a partial fraction with all even entries. However, if $p$ and $q$ have opposite parity, then it is easy to see that the remainders $r_0, r_1, \dots, r_n$ alternate in parity, and because $p$ and $q$ have opposite parity their greatest common divisor, $r_n$, is odd. This means that $a_n$ is even. If $p$ and $q$ are both odd, then choosing $a_1$ to be odd will cause $r_2$ to be even and $r_1$ and $r_2$ are now of opposite parity. Thus continuing to choose even partial quotients from this point on will end with the last partial quotient being even.

To carry out the Euclidean algorithm, we do not require that $p$ and $q$  be relatively prime. But when considering a fraction that represents a 2-bridge knot or link, we assume the fraction is reduced. Thus we will never have the case that both $p$ and $q$ are even. We have proven the following lemma, except for the final assertion, which is easily verified by induction.

\begin{lemma}
\label{special fraction} Let $\frac{p}{q}$ be a reduced fraction. Then we may express $p/q$ as 
$$\frac{p}{q}=r+[a_1, a_2, \dots, a_n],$$
where each $a_i$ is a nonzero even integer. 
If $q$ is odd, this form is unique.  Moreover, $n$ must be even and $p$ and $r$ have the same parity. 
If $q$ is even then there are exactly two ways to express $p/q$ in this way. In one, $r=\lfloor \frac{p}{q} \rfloor$ and in the other $r=\lceil \frac{p}{q} \rceil$. The number $n$ of partial quotients is not necessarily the same in each case, but is always odd.
\end{lemma}

The Euclidean algoritihm described above will never produce $a_i=0$ since $|r_{i-1}|>|r_{i+1}|$. However, we can easily make sense of continued fractions that use zeroes. In this case it is not difficult to show that a zero can be introduced or deleted from a continued fraction as follows:
$$r+[\dots, a_{k-2}, a_{k-1}, 0, a_{k+1}, a_{k+2}, \dots]=r+[\dots, a_{k-2}, a_{k-1}+a_{k+1}, a_{k+2},\dots].$$
Using this property, every continued fraction with all even partial quotients can be expanded so that each partial quotient is either $-2, 0$, or $2$. For example, a partial quotient of $6$ would be expanded to $2,0,2,0,2$ and $-4$ to $-2,0,-2$. This leads us to the following definition.

\vspace{-.25in}

\begin{definition} {\rm Let $\cal S$ be the set of all integer vectors $(a_1, a_2, \dots, a_n)$ such that

\vspace{-.25in}

\begin{enumerate}
\item each $a_i\in \{-2, 0, 2\}$, 
\item $a_1\ne0$ and $a_n\ne 0$,
\item if $a_i=0$ then $a_{i-1}=a_{i+1}=\pm 2$.
\end{enumerate}

\vspace{-.25in}

We call $\cal S$ the set of {\it expanded even} vectors.}
\end{definition}

\vspace{-.25in}

We may further define an equivalence relation on $\cal S$ by declaring that ${\bf a}, {\bf b} \in \cal S$ are equivalent if ${\bf a}=\pm {\bf b}$ or ${\bf a}=\pm {\bf b}^{-1}$ where ${\bf b}^{-1}$ is ${\bf b}$ read backwards.  We denote the equivalence class of ${\bf a}$ as $\hat{\bf a}$ and the set of all equivalence classes as $\hat{ \cal S}$. 
Notice that $\hat{\cal S}=\hat{\cal S}_{even} \cup \hat{\cal S}_{odd}$, 
where $\hat{\cal S}_{even}$ consists of classes represented by vectors of even length and $\hat{\cal S}_{odd}$ consists of classes represented by vectors of odd length.

Because every fraction has a (nearly) canonical representation as a continued fraction using a vector of even partial quotients, we can prove the following proposition. This provides a nice characterization of 2-bridge knots and links which almost certainly has already appeared somewhere in the literature.  We will see in the next sections that working with expanded even vectors (as opposed to simply vectors of nonzero even integers, or some other variant of the Euclidean algorithm) is a natural choice when working with the Ohstuki-Riley-Sakuma partial order.

\begin{proposition}
\label{2-bridge classification}
Let $\Phi: \hat{\cal S} \to \cal B$ be given by $\Phi(\hat{\bf a})=K_{p/q}$ where $p/q=0+[{\bf a}]$. Then 

\vspace{-.25in}

\begin{enumerate}
\item The restriction of $\Phi$ to $\hat{\cal S}_{even}$ is a bijection onto the set of equivalence classes of 2-bridge knots.
\item The restriction of $\Phi$ to $\hat{\cal S}_{odd}$ is a two-to-one map onto the set of equivalence classes of 2-bridge links.
\end{enumerate}
\end{proposition}

\noindent{\bf Proof:} Suppose that ${\bf a}=(a_1, a_2, \dots, a_n)$. It is well known, and not hard to prove, that if $p/q=0+[\bf a]$, then $0+[{\bf -a}]=-p/q$ and $0+[{\bf a}^{-1}]=p'/q$ where $p p' \equiv (-1)^{n+1} (\mbox{mod } q)$. Thus $\Phi$ is well defined.

Suppose now that $\hat {\bf a}, \hat{\bf b} \in \hat{\cal S}_{even}$ and that $\Phi(\hat {\bf a})=\Phi(\hat{\bf b})$. If $p/q=0+[{\bf a}]$ and $p'/q'=0+[{\bf b}]$, then $q=q'$ and either $p\equiv \pm p' \ (\mbox{mod } q)$, or $p p'\equiv \pm 1 \ (\mbox{mod } q)$. 

If $p\equiv \pm p' \ (\mbox{mod } q)$, then $\pm p'/q=p/q+r$ for some integer $r$ and the fraction $\pm p'/q$ can be expressed as both $0+[\pm {\bf b}]$ and $r+[{\bf a}]$. At this point we would like to apply  Lemma~\ref{special fraction}, but are faced with the technical difficulty that $\bf a$ and $\bf b$ are expanded even vectors rather than just vectors of nonzero even integers. Hence we must first {\it contract}, that is ``unexpand,''  both $\bf a$ and $\bf b$ (this process is unique) to nonzero vectors $\bf a'$ and $\bf b'$, respectively. It is still the case that $\pm p'/q$ can be expressed as both $0+[\pm {\bf b'}]$ and $r+[{\bf a'}]$ and now by Lemma~\ref{special fraction}, such an expression is unique. So we have $r=0$ and $\pm {\bf b'}={\bf a'}$. Hence $\pm {\bf b}={\bf a}$ and $\hat{\bf a}=\hat{\bf b}$.

If instead, $p p'\equiv \pm 1 \ (\mbox{mod } q)$, consider $p''/q=0+[{\bf b}^{-1}]$. Because the length of $\bf b$ is even,  $p' p'' \equiv -1 \ (\mbox{mod } q)$ and hence $p p'\equiv \mp p' p'' \ (\mbox{mod } q)$. Since $p'$ and $q$ are relatively prime, we may cancel $p'$ and obtain $ p \equiv \mp p''\ (\mbox{mod } q)$. Proceeding as in the first case, we now have $\pm {\bf b}^{-1}={\bf a}$ and $\hat{\bf a}=\hat{\bf b}$.

A similar argument can be given in the case of links.
\hfill $\square$

\section{The Ohtsuki-Riley-Sakuma partial order}
\label{ORS partial order section}

In \cite{ORS}, Ohtsuki, Riley and Sakuma systematically construct epimorphisms between 2-bridge knots and links which preserve peripheral structure. 
In particular, they show that if  $J=\Phi(\hat{\bf c})$ where the vector {\bf c} has the form
$${\bf c}=(\epsilon_1 {\bf a}, {\bf 2 c_1}, \epsilon_2 {\bf a}^{-1}, {\bf 2  c_2}, \epsilon_3 {\bf a}, {\bf 2  c_3}, \dots, \epsilon_n  {\bf a}^{(-1)^{n-1}}),$$
with $\epsilon_i \in \{-1, 1\}$, $c_i$  an integer, and  furthermore, if $c_i=0$ then $\epsilon_i=\epsilon_{i+1}$, then there is an epimorphism from the group of $J$ to the group of $K=\Phi(\hat{\bf a})$. Here we have used $({\bf a}, {\bf b})$ to represent the concatenation of {\bf a} and {\bf b} and we denote by ${\bf 2  c_i}$ the vector $\pm (2, 0, 2, 0, \dots, 2)$ whose entries add to $2 c_i$. Ohtsuki, Riley and Sakuma's result motivates the following definition.

\begin{definition} {\rm Let ${\bf a}, {\bf c} \in  \cal S$. We say that $\bf c$ {\it admits a parsing with respect to} $\bf a$ if 
$${\bf c}=(\epsilon_1{\bf a}, {\bf2  c_1}, \epsilon_2 {\bf a}^{-1}, {\bf 2  c_2}, \epsilon_3 {\bf a}, {\bf 2  c_3}, \dots, \epsilon_n  {\bf a}^{(-1)^{n-1}}),$$
where $\epsilon_1=1$, $\epsilon_i \in \{-1, 1\}$ for $i>1$, and each $c_i$ is an integer. Furthermore, if $c_i=0$ then $\epsilon_i=\epsilon_{i+1}$. We say that $\hat {\bf c}$ admits a {\it parsing with respect to $\hat{\bf a}$} if some representative of $\hat {\bf c}$ admits a parsing with respect to some representative of $\hat{\bf a}$.
} 
\end{definition}
Note that if {\bf c} parses with respect to {\bf a}, then this parsing is unique. This can easily be proven by induction on the length of {\bf c}. 

If {\bf c} parses with respect to {\bf a}, then we call the vectors ${\bf 2  c_i}$ {\bf a}-{\it connectors}. The {\bf a}-connectors separate the {\it {\bf a}-tiles} $\epsilon_1 {\bf a},  \epsilon_2 {\bf a}^{-1}, \dots,$ and $ \epsilon_n  {\bf a}^{(-1)^{n-1}}$. The set of all possible connectors is $${\cal C}=\{ (0), \pm(2), \pm(2,0,2), \pm(2,0,2,0,2), \dots \}.$$ These vectors will be important not only in their role as connectors in a given parsing, but more generally as the {\it components} from which all vectors in $\cal S$ can be built (by concatenation). Given any vector ${\bf c} \in {\cal S}$, and a substring {\bf C} contained in {\bf c}, we say that {\bf C} is a {\it maximal component in {\bf c}} if it is an element of $\cal C$  and is not contained in any larger substring of {\bf c} which is an element of $\cal C$. Notice that ${\bf C}=(0)$ is never a maximal component in {\bf c}, unless ${\bf c}=(0)$. Clearly every vector ${\bf c}\in {\cal C}$ can be uniquely decomposed into a sequence of maximal components. 

\begin{definition} \label{ORS order def}If $K_1=\Phi(\hat{\bf c})$ and $K_1=\Phi(\hat{\bf a})$ are 2-bridge knots, we define $K_1\ge K_2$ if  $\hat{\bf c}$ parses with respect to $\hat{\bf a}$. 
\end{definition}

It is not hard to show that Definition~\ref{ORS order def} does indeed define a partial order on the set of 2-bridge knots. We call this the {\it Ohtsuki-Riley-Sakuma partial order}. Notice that if $K_1\ge K_2$, then there exist an epimorphism which preserves peripheral structure from the group of $K_1$ onto the group of $K_2$, and hence, $K_1$ is greater than $K_2$ in the Silver-Whitten partial order. Thus the Silver-Whitten partial order is a refinement of the Ohtsuki-Riley-Sakuma partial order.

We say that $J$ is an {\it upper bound} of of a set of 2-bridge knots $\{K_i\}_{i \in I}$ if $J\ge K_i$ for all $i\in I$. The knot $J$ is a {\it least} upper bound if it is an upper bound and not strictly greater than any other upper bound of  $\{K_i\}_{i \in I}$. Note that it is possible for a set of knots to have no least upper bounds, one least upper bound, or multiple least upper bounds. Note also that if  $\{K_i\}_{i \in I}$ has an upper bound, then $I$ is finite. This follows from Simon's Conjecture, which was recently proven in \cite{AL:2010}, but which was previously known in the special case of 2-bridge knots \cite{BBRW:2009}.

We now consider an example that illustrates a key  idea of this paper. 

{\bf Example 1:} Consider the 2-bridge knots $K_1=K_{4/7}, K_2=K_{24/41}$ and $J=K_{322892/551327}$.  We have  $K_1=\Phi(\hat {\bf a})$ with ${\bf a}=(2,-2,0,-2)$,  $K_2=\Phi(\hat{\bf b})$ with ${\bf b}=(2,-2,0,-2,2,-2,0,-2)$ and $J=\Phi(\hat{\bf c})$ with 
$${\bf c}=(2,-2,0,-2,2,-2,0,-2,2,2,0,2,-2,2,0,2,-2,-2,0,-2,2,-2,0,-2,2,-2,0,-2).$$ 
Since we can parse $\bf c$ as
$$\overbrace{[2,-2,0,-2]}^{{\bf a}},(2),\overbrace{[-2,0,-2,2]}^{{\bf a}^{-1}},(2,0,2),\overbrace{[-2,2,0,2]}^{{-\bf a}},(-2),\overbrace{[-2,0,-2,2]}^{{\bf a}^{-1}},(-2,0,-2),\overbrace{[2,-2,0,-2]}^{{\bf a}}$$ as well as 
$$\overbrace{[2,-2,0,-2,2,-2,0,-2]}^{{\bf b}},(2),\overbrace{[2,0,2,-2,2,0,2,-2]}^{{-\bf b}^{-1}},(-2,0,-2),\overbrace{[2,-2,0,-2,2,-2,0,-2]}^{{\bf b}}$$
it follows that $J$ is an upper bound for the set $\{K_1, K_2\}$. 
In the first parsing we see five copies of the tile $\bf a$ laid down alternately forward and backward, separated by the {\bf a}-connectors 
$(2), (2,0,2), -(2)$, and $-(2,0,2)$.  To emphasize the difference between tiles and connectors, we have used square and round brackets, respectively. Note that the third {\bf a}-tile has been negated.  
In the second parsing we see three copies of the tile $\bf b$ separated by two {\bf b}-connectors, and with the second {\bf b}-tile negated.  As this example illustrates, whenever a vector {\bf c} can be parsed in two or more different ways, it represents an upper bound for some set of distinct knots.

\begin{figure}
    \begin{center}
    \leavevmode
    \scalebox{.50}{\includegraphics{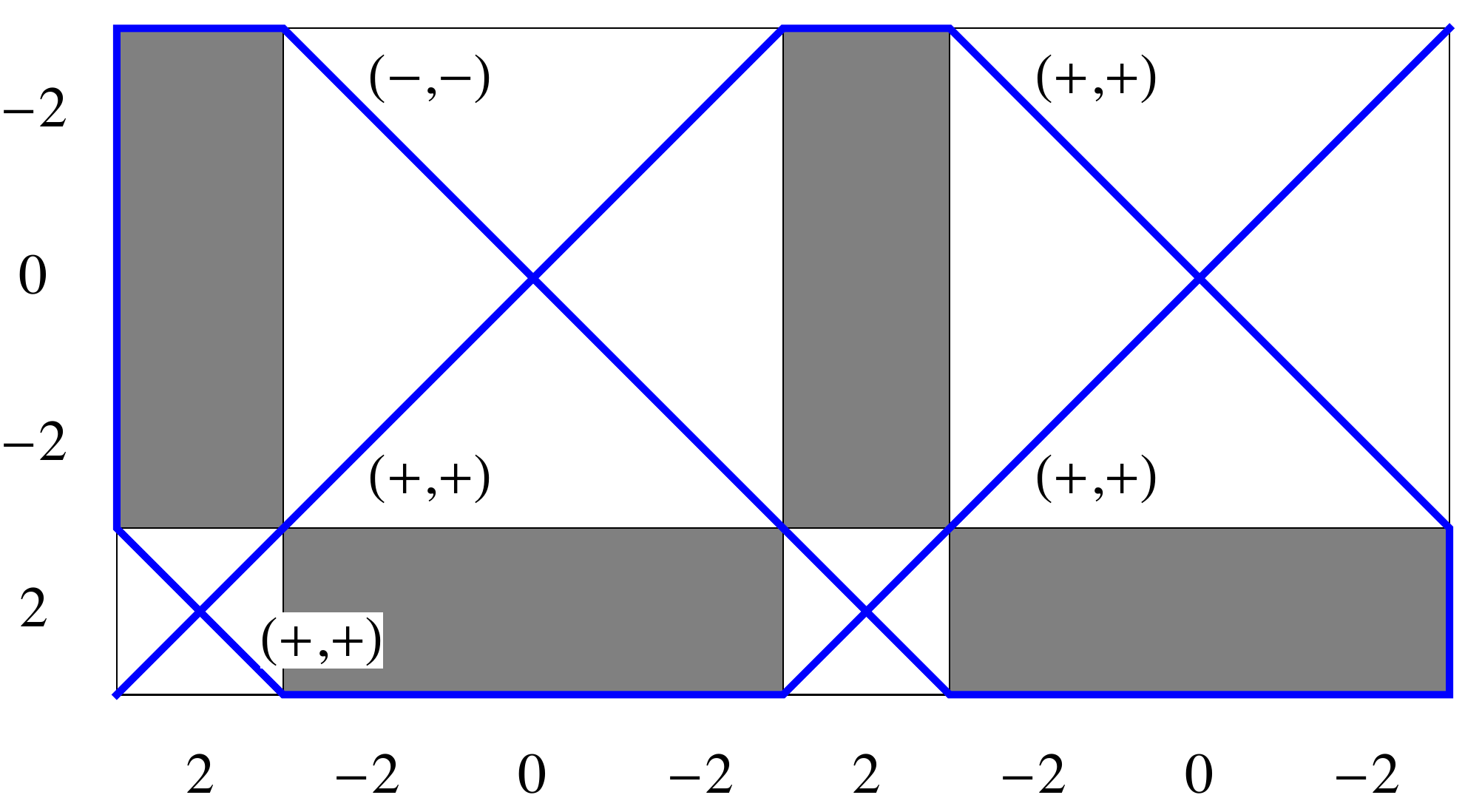}}
    \end{center}
\caption{The upper bound {\bf c} of $\{{\bf a}, {\bf b}\}$, seen as a path in ${\bf a}\times {\bf b}$. }
\label{cPath1}
\end{figure}

There is a nice geometric interpretation of the vectors {\bf a}, {\bf b} and {\bf c} in Example 1 which we display in Figure~\ref{cPath1}. Here we form the {\it product} of {\bf a} and {\bf b}, with {\bf a} labeling the rows of the product, from bottom to top, and {\bf b} the columns, from left to right.  The vector {\bf c} is now represented by a path that starts at the lower left corner and moves  along lines of slope $1$ or $-1$, ending at the upper right corner. When the path hits an edge of the product, it moves along the edge, traversing a connector. The {\bf a}-connectors appear on the upper and lower edges; the {\bf b}-connectors appear on the right and left edges. Additionally, each diagonal segment of the path {\bf c} is labeled with a pair of signs, $(\epsilon, \eta)$, as follows. The first sign indicates that the elements of {\bf c} along that diagonal lie in $\epsilon {\bf a}^{\pm1}$ in the parsing of {\bf c} with respect to {\bf a}. Similarly,  the second sign indicates that these elements of {\bf c} lie in $\eta {\bf b}^{\pm 1}$. If the row and column labels of a cell containing a diagonal segment of {\bf c} are $x$ and $y$, then $x=\epsilon \eta y$.  Hence two diagonal segments of {\bf c} that intersect must be labeled with signs whose product is the same. In this example,  the diagonal segments of {\bf c} form a connected subset of $\mathbb R^2$. Since the initial segment is labeled $(+,+)$, we see that every segment must be labeled either $(+,+)$ or $(-,-)$. As we traverse the path {\bf c}, if two consecutive diagonal segments both lie in the same {\bf a} or {\bf b}-tile, the corresponding sign cannot change, and hence the pair of signs cannot change. Thus a change from $(+,+)$ to $(-,-)$ or back can only occur at a point where the path {\bf c} passes through a corner of the product. In the example shown in Figure~\ref{cPath1}, the path traverses the product a total of five times in the vertical direction and three times in the horizontal direction since there are five {\bf a}-tiles and three {\bf b}-tiles in the two parsings. It is not the case that every vector {\bf c} that can be parsed with respect to both {\bf a} and {\bf b} can be depicted as a path in ${\bf a} \times {\bf b}$ in this way. But we shall prove in Section~\ref{ORS partial order} that this is the case if {\bf c} is minimal in length. Finally, we have shaded the product checkerboard fashion with respect to the maximal component decompositions of {\bf a} and {\bf b}.

At this point we can see the value of using expanded even vectors instead of nonzero even vectors. If a connector is zero, and is removed by contracting, the parsing can be hidden. For example, if
$${\bf c}=(2,4,4,4,4,4,4,4,4,4,4,4,4,4,4,2),$$
then  it is not immediately apparent that $\Phi(\hat {\bf c})$ is an upper bound of both $\Phi(\hat{\bf a})$ and $\Phi(\hat {\bf b})$ where ${\bf a}=(2,4,4,2)$ and ${\bf b}=(2,4,4,4,4,2)$. However, after expanding ${\bf a}, {\bf b}$ and $\bf c$, a double parsing like the previous example will reveal itself.

 A key observation of this paper is the following theorem,  which says that every pair of 2-bridge knots with an upper bound arises as in Example 1, by parsing a vector in ${\cal S}_{even}$ in two different ways.

\begin{theorem} \label{double parsing theorem}The 2-bridge knot  $K$ is an upper bound of both $K_1$ and $K_2$ with respect to the  Ohtsuki-Riley-Sakuma partial order if and only if there exists $\hat{\bf a}, \hat{\bf b},  \hat{\bf c} \in \hat{\cal S}_{even}$ such that $K_1=\Phi(\hat{\bf a}),  K_2=\Phi(\hat{\bf b}), K=\Phi(\hat{\bf c})$, and $\hat{\bf c}$ can be parsed with respect to both $\hat{\bf a}$ and $\hat{\bf b}$.
\end{theorem}

\noindent{\bf Proof:} Suppose $K_1$ corresponds to the vector ${\bf a} \in  {\cal S}_{even}$. If $K\ge K_1$, then it corresponds to a vector ${\bf c} \in {\cal S}_{even}$ of the form  ${\bf c}=(\epsilon_1 {\bf a}, {\bf 2  c_1}, \epsilon_2 {\bf a}^{-1}, {\bf 2  c_2}, \epsilon_3 {\bf a}, {\bf 2  c_3}, \dots, \epsilon_{2n+1}  {\bf a})$.
Similarly, suppose $K_2$ corresponds to the vector ${\bf b}\in {\cal S}_{even}$. Then $K$ must also correspond to a vector $\bf d$ of the form 
${\bf d}=(\delta_1 {\bf b}, {\bf 2  d_1}, \delta_2 {\bf b}^{-1},  {\bf 2  d_2}, \delta_3 {\bf b}, {\bf 2  d_3}, \dots, \delta_{2m+1}  {\bf b)}$. Thus by Proposition~\ref{2-bridge classification}, we must have $\hat {\bf c}=\hat{\bf d}$ and $\bf c$ parses with respect to both $\hat{\bf a}$ and $\hat{\bf b}$.
\hfill $\square$

In Section~\ref{ORS partial order} we will completely characterize those  pairs of knots that have an upper bound. However, we state the following simple corollary of Theorem~\ref{double parsing theorem} now  because it settles the  question that originally motivated this research. 

\begin{corollary} \label{knots of same length} Suppose $K_1=\Phi(\hat{\bf a})$ and $K_2=\Phi(\hat{\bf b})$ are knots. If $K_1$ and $K_2$  have an upper bound with respect to the Ohtsuki-Riley-Sakuma partial order and  $\hat {\bf a}$ and $\hat {\bf b}$ have the same length, then  $\hat {\bf a}=\hat {\bf b}$. 
\end{corollary}

\noindent{\bf Example 2:} By Corollary \ref{knots of same length}, the trefoil knot $3_1$, which is given by $(2,-2)$, and the figure eight knot $4_1$, which is given by $(2,2)$, do not have an upper bound.

Because the map $\Phi$ is two-to-one onto the set of 2-bridge links, the situation is more complicated for links. For example, consider  the 2-bridge link $K_{11/30}=\Phi(\hat{\bf c})=\Phi(\hat{\bf d})$ where ${\bf c}=(2,2,-2,2,2)$ and ${\bf d}=(2,-2,-2,-2,2)$. Clearly $\bf c$ parses with respect to $(2,2)$ and $\bf d$ parses with respect to $(2,-2)$. Hence there is an epimorphism of the group of the link $K_{11/30}$ onto  both the trefoil and the figure eight knots. However, neither $\bf c$ nor $\bf d$ can be parsed with respect to {\sl both} $(2,2)$ and $(2,-2)$.  We note that $K_{11/30}$ is one member of the family $K_\frac{9+20 n}{20+50 n}$, all of which have groups that map onto both  the trefoil and the figure eight groups, provided $n \equiv -1 \ (\mbox{mod } 3)$.

\pagebreak
\section{Comparability, lower, and upper bounds}\label{ORS partial order}

Given two 2-bridge knots, $K_1$ and $K_2$, it is straightforward to check if they are comparable in the Ohtsuki-Riley-Sakuma  partial order. In order to decide if $K_1\ge K_2$, we first represent each as vectors {\bf a} and {\bf b}, respectively, in ${\cal S}_{even}$  and then check to see if $\hat{\bf a}$ parses with respect to $\hat{\bf b}$. In fact, given a vector ${\bf a}\in {\cal S}_{even}$, it is easy to find all vectors {\bf c} such that {\bf a}  parses with respect  to {\bf c} and, thus, to determine all 2-bridge knots $K$ such that $K_1 \ge K$. 

It is also an easy matter to find all lower bounds for a given set of knots. Of course every 2-bridge knot is greater than the unknot, so the unknot is a lower bound for any set of 2-bridge knots. But given a finite set $\{K_i\}_{i=1}^n$ of 2-bridge knots, we can first determine all knots smaller than each $K_i$. The intersections of these sets then provides the complete set of lower bounds. Determining which lower bounds are greatest lower bounds is then straightforward since we need only check a finite number of pairs of knots for comparability.

The problem of finding upper bounds for a set of knots is much more difficult. The rest of this section will be devoted to proving the following theorem.

\begin{theorem}\label{main theorem} The 2-bridge knots $K_1=\Phi(\hat{\bf a})$ and $K_2=\Phi(\hat{\bf b})$ are incomparable and $\{K_1, K_2\}$ has an upper bound with respect to the Ohtsuki-Riley-Sakuma partial order if and only if
$${\bf a}=(\overbrace{{\bf w}, {\bf w}, \dots, {\bf w}}^{p}, {\bf e}) \quad \mbox{and} \quad {\bf b}=(\overbrace{{\bf w}, {\bf w}, \dots, {\bf w}}^q, {\bf e}),$$
where $\bf e$ is some (possibly empty) vector in ${\cal S}_{even}$,  ${\bf w}=({\bf e}, {\bf m}, {\bf e}^{-1}, {\bf n})$, $m$ and $n$ are even integers, and finally, neither $2q+1$ nor $2p+1$ divides the other.
\end{theorem}

Notice that if we denote the vector $(\overbrace{{\bf w}, {\bf w}, \dots, {\bf w}}^{p}, {\bf e})$ as $({\bf w}^p, {\bf e})$, then it is not hard to show that $({\bf w}^p, {\bf e})$ parses with respect to $({\bf w}^q, {\bf e})$ if and only if $2q+1$ divides $2p+1$. This observation makes one direction of the theorem  straightforward. Suppose ${\bf a}=({\bf w}^{p}, {\bf e})$, 
${\bf b}=({\bf w}^{q}, {\bf e})$,  $\bf e$ is a (possibly empty) vector in 
${\cal S}_{even}$,  and ${\bf w}=({\bf e}, {\bf m}, {\bf e}^{-1}, {\bf n})$ where $m$ and $n$ are even integers. Because both $2p+1$ and $2q+1$ divide $2(2 pq+p+q)+1$, we have that $({\bf w}^{2 pq+p+q}, {\bf e})$ parses with respect to both $({\bf e}^p, {\bf e})$  and $({\bf w}^q, {\bf e})$. Hence $\Phi(({\bf w}^{2 pq+p+q}, {\bf e}))$ is an upper bound of $\{K_1, K_2\}$.

We will prove the other direction of Theorem~\ref{main theorem} by first proving the following proposition which is, in some sense, a special case of the theorem.  We then state and prove a number of technical lemmas that allow us to reduce the general case to that of Proposition~\ref{gcd result}. 

\pagebreak
\begin{proposition}\label{gcd result} Suppose ${\bf a}, {\bf b}, {\bf c} \in {\cal S}_{even}$, {\bf c} parses with respect to both {\bf a} and {\bf b}, {\bf a} does not parse with respect to $\hat{\bf b}$, and {\bf b} does not parse with respect to $\hat{\bf a}$. If ${\bf c}=(C_1, C_2, \dots, C_{t-1})$ where each $C_i \in {\cal C}$  and there exist $r$ and $s$ such that
\vspace{-.25 in}
\begin{enumerate}
\item $t=\mbox{lcm}(r, s)$, 
\item ${\bf a}=(C_1, C_2, \dots, C_{r-1})$,
\item ${\bf b}=(C_1, C_2, \dots, C_{s-1})$, 
\item the parsing of {\bf c} with respect to {\bf a} is given by

$\begin{array}{l}
[C_1, C_2, \dots, C_{r-1}]\, (C_r) \, [C_{r+1}, C_{r+2}, \dots, C_{2r-1}]\, (C_{2r})\, \dots [C_{t-r+1}, C_{t-r+2}, \dots, C_{t-1}]\\
=(\epsilon_1 {\bf a}, C_r, \epsilon_2 {\bf a}^{-1}, C_{2 r}, \epsilon_3 {\bf a}, C_{3r}, \dots, C_{t-r}, \epsilon_\frac{s}{d}{\bf a})
\end{array}$

where each $\epsilon_i =\pm 1$ and $\epsilon_1=1$, 

\item the parsing of {\bf c} with respect to {\bf b} is given by

$\begin{array}{l}
[C_1, C_2, \dots, C_{s-1}]\, (C_s) \, [C_{s+1}, C_{s+2}, \dots, C_{2s-1}]\, (C_{2s})\, \dots [C_{t-s+1}, C_{t-s+2}, \dots, C_{t-1}]\\
=(\eta_1{\bf b}, C_s, \eta_2 {\bf b}^{-1}, C_{2 s}, \eta_3 {\bf b}, C_{3s}, \dots, C_{t-s}, \eta_\frac{r}{d}{\bf b})
\end{array}$

where each $\eta_i =\pm 1$ and $\eta_1=1$, and 
\item whenever  $C_k$ lies in both the {\bf a}-tile $\epsilon_i {\bf a}^{\pm 1}$ and  the {\bf b}-tile $\eta_j {\bf b}^{\pm 1}$, then $\epsilon_i=\eta_j$,

\end{enumerate}
\vspace{-.25 in}
then ${\bf a}=({\bf w}^{p}, {\bf e})$ and 
${\bf b}=({\bf w}^{q}, {\bf e})$,  where $\bf e$ is a (possibly empty) vector in 
${\cal S}_{even}$,  ${\bf w}=({\bf e}, {\bf m}, {\bf e}^{-1}, {\bf n})$, $m$ and $n$ are even integers, and $p$ and $q$ are natural numbers. Furthermore, neither $2p+1$ nor $2q+1$ divides the other.
\end{proposition} 
{\bf Proof:} We cannot have $r=s$ since neither of {\bf a} nor {\bf b} parses with respect to the other. Without loss of generality, assume that $r<s$. Let $d=\gcd(r,s)$. If $d=r$ then {\bf b} parses with respect to {\bf a}.  Hence, $0<d<r<s<t$. Since {\bf a, b} and {\bf c} are in ${\cal S}_{even}$, it follows that $r, s, t, d, r/d$ and $s/d$ are all odd.

{\bf Claim:} If $0<i<r$, $0<j<s$ and $j=i+2 u d$ for some positive integer $u$, then $C_i=C_j$.

To prove this claim, form the vector ${\bf c}'=({\bf c}, 0, {\bf c}^{-1})=(C'_1, C'_2, \dots, C'_{2t-1})$ with 
$$C'_k=\left\{ \begin{array}{lll}C_k &\mbox{if}&0<k<t, \\
0&\mbox{if}& k=t,\\
C_{2t-k}&\mbox{if}& t<k<2t. \end{array}\right.$$
Note that ${\bf c}'$ parses with respect to both {\bf a} and {\bf b}. Furthermore, if $k$ is not a multiple of $r$, then $C'_k$ lies in an {\bf a}-tile, and the sign of that {\bf a}-tile is $\epsilon_{\sigma(k)}$ where 
$$\sigma(k)=\left\{ \begin{array}{lll}\lceil \frac{k}{r} \rceil&\mbox{if}&0<k<t,\\
\\
\lceil \frac{2t-k}{r} \rceil&\mbox{if}&t<k<2t.
  \end{array}\right.$$
  Similarly, if $j$ is not a multiple of $s$, then $C'_j$ lies in a {\bf b}-tile, and the sign of that {\bf b}-tile is $\eta_{\tau(j)}$ where 
$$\tau(j)=\left\{ \begin{array}{lll}\lceil \frac{j}{s} \rceil&\mbox{if}&0<j<t,\\
\\
\lceil \frac{2t-j}{s} \rceil&\mbox{if}&t<j<2t.
  \end{array}\right.$$
  If neither $r$ nor $s$ divide $n$, then $C'_n$ lies in both an {\bf a}-tile and a {\bf b}-tile and part six of our hypothesis gives $\epsilon_{\sigma(n)}=\eta_{\tau(n)}$.
  
Using the {\bf a}-parsing of ${\bf c}'$, we see that  if $0<i<r$ and $0<k<s/d$, then $$C'_{i+2kr}=\epsilon_{\sigma(i+2kr)}C_i.$$  Similarly, using the {\bf b}-parsing, we see that if $0<j<s$ and $0<l<r/d$, then
$$C'_{j+2ls}=\eta_{\tau(j+2ls)}C_j.$$

We want to show that there exist $k$ and $l$ with $0\le k <\frac{s}{d}$ and $0\le l<\frac{r}{d}$ such that $$i+2 kr=j+2 l s$$ or equivalently
$$kr-ls=ud.$$ This will imply that $$C_i=\epsilon_{\sigma(i+2kr)}C'_{i+2kr}=\eta_{\tau(j+2ls)}C'_{j+2ls}= C_j.$$
Now since $\gcd(r, s)=d$ there certainly exist integers $x$ and $y$ such that $xr-ys=ud$. It is easy to show that the set of all such integers $x$ and $y$ are given by $$(x,y)=(x_0,y_0)+v \left (\frac{s}{d}, \frac{r}{d}\right),$$ where $(x_0, y_0)$ is one such pair and $v$ is any integer. Thus, we may select $v_0$ so that $k=x_0+v_0 \frac{s}{d}$ is contained in the half open interval $[0, \frac{s}{d})$. Let $l=y_0+v_0 \frac{r}{d}$. Since $kr-ls=ud$, and $ud>0$, we have $ls<kr$. Hence, $ls<kr<\frac{s}{d} r$ and so $l<\frac{r}{d}$. We also know that $0\le l$, for if not, and $l<0$ then in fact $l\le -1$. But this implies $ud=kr-ls>-ls\ge s$, which is impossible since $j=i+2 ud<s.$ This completes the proof of the claim.

Now let ${\bf e}=(C_1, C_2, \dots, C_{d-1}), {\bf m}=C_d, {\bf E}=(C_{d+1}, C_{d+2}, \dots, C_{2d-1})$ and ${\bf n}=C_{2d}$. We have that 
$${\bf a}=(({\bf e, m, E, n})^p,{\bf e})$$ and $${\bf b}=(({\bf e, m, E, n})^{q},{\bf e}),$$ where $p={(r/d-1)/2}$ and $q=(s/d-1)/2$.
We are assuming that {\bf b} is longer than {\bf a}. Hence the first {\bf a}-connector, $C_r={\bf m}$, lies in the first {\bf b}-tile and is followed in the {\bf b}-tile by {\bf E}. However, it is followed in the {\bf a}-parsing by ${\bf e}^{-1}$. This implies that ${\bf E}={\bf e}^{-1}$ and thus {\bf a} and {\bf b} have the desired form.  Finally, $2p+1=r/d$ and $2q+1=s/d$. Since $r/d$ and $s/d$ are relatively prime, and each is at least 3, neither divides the other.  \hfill $\square$

Our goal now is to reduce Theorem~\ref{main theorem} to the case of Proposition~\ref{gcd result}. 
Suppose that $K_1=\Phi(\hat{\bf a})$ and $K_2=\Phi(\hat{\bf b})$ are incomparable knots and $K=\Phi({\hat{\bf c}})$ is an upper bound of $\{K_1, K_2\}$. Then there are representatives of $\hat{\bf a}, \hat{\bf b}$, and $\hat{\bf c}$, which we call {\bf a}, {\bf b}, and {\bf c}, respectively, such that {\bf c} parses with respect to both {\bf a} and {\bf b}. Furthermore, we may assume that {\bf c} is a shortest representative among all  upper bounds of $\{K_1, K_2\}$ and that {\bf a} is shorter than {\bf b}.  We need to show that {\bf c} can be decomposed into a sequence of components satisfying the hypotheses of Proposition~\ref{gcd result}. This will require several technical lemmas. 

Before stating our first lemma, which will establish several properties of the vector {\bf c}, we introduce some needed terms. If $k$ is less than the length of $\bf c$ and the $k$-th entry of $\bf c$ marks the end of either a connector or a tile in both the $\bf a$ and {\bf b}-parsings,  we say that there is a {\it seam} between elements $k$ and $k+1$ of $\bf c$. If, in both parsings, tiles are on the same side of the seam and connectors on the other, we call the seam a {\it pure} seam, otherwise it is a {\it mixed} seam. Note that at a mixed seam, neither of the two connectors can be $(0)$ since this would imply that one of the tiles begins or ends in a zero.  
In Example 1,  mixed seams occur after the 9th and 17th entries of $\bf c$. In the geometric interpretation of {\bf c} as a path in ${\bf a} \times {\bf b}$, mixed seams correspond to the path traveling around the corner of the product. As mentioned before, mixed seams represent the only opportunities to negate a tile or its inverse.

\begin{lemma} \label{properties of c}Suppose ${\bf a}, {\bf b}, {\bf c} \in {\cal S}_{even}$,  {\bf c} parses with respect to both ${\bf a}$ and ${\bf b}$, and that no shorter vector in ${\cal S}_{even}$ has this property. Then all of the following are true.
\vspace{-.25in}
\begin{enumerate}
\item The vector {\bf c} contains no pure seams.
\item The intersection between any {\bf a} and {\bf b}-connector is empty. 
\item  Every {\bf a}-connector is either a maximal component or the central zero of a maximal component in a {\bf b}-tile.
\item Every {\bf b}-connector is either a maximal component or the central zero of a maximal component in an {\bf a}-tile. 
\end{enumerate}
\end{lemma} 
{\bf Proof:} First note that if an {\bf a}-connector and a {\bf b}-connector intersect and are not equal, then their right or left ends do not align. This will cause either an {\bf a}-tile or a {\bf b}-tile to begin or end with a zero unless one connector is $(0)$ and is centered in the other. Now suppose that {\bf c} contains a pure seam. The {\bf a} and {\bf b}-connectors adjacent to the seam intersect and are aligned at one end, hence must be equal. Thus pure seams exist on both sides of these coincident connectors. Let {\bf x} and {\bf y} be the substrings of {\bf c} which lie before the first seam and after the second seam, respectively. Since {\bf a}, {\bf b}, and {\bf c} are all of even length, then  exactly one of {\bf x} or {\bf y} is in $S_{even}$. If ${\bf x} \in {\cal S}_{even}$, then it parses with respect to both {\bf a} and {\bf b} and is shorter than {\bf c}, a contradiction. If instead, ${\bf y}\in S_{even}$, then {\bf y} or $-{\bf y}$, parses with respect to both ${\bf a}$ and ${\bf b}$ and is shorter in length. (Note that {\bf y} cannot  begin with $-{\bf a}$ and ${\bf b}$ for example. If it  did, then  by considering {\bf y}, we see that the first element of {\bf a}  is the negative of the first element of {\bf b}. But by considering {\bf x}, we see that the first element of {\bf a} equals the first element of {\bf b}. Thus the first elements of {\bf a} and {\bf b} must be zeroes, a contradiction.) This establishes the first property.

Suppose now that an {\bf a} and {\bf b}-connector intersect. Since there are no pure seams, they cannot be equal and so one must be zero and centered in the other. Without loss of generality, assume there exists an {\bf a}-connector ${\bf d}=(0)$ centered in a {\bf b}-connector ${\bf d}'=\pm(2, \dots, 0 \dots, 2)$. Let {\bf s} be the largest substring of {\bf c} centered at {\bf d} that parses with respect to $\hat{\bf a}$ and such that $\mbox{abs}({\bf s})=\mbox{abs}({\bf s}^{-1})$, where $\mbox{abs}({\bf s})$ is the vector obtained by negating all negative entries of {\bf s}. Let {\bf t} be similarly  defined with respect to the {\bf b}-parsing. Let $\sigma_{\bf s}$ and $\sigma_{\bf t}$ be the involutions of {\bf s} and {\bf t},  respectively, centered at {\bf d}, that take each of these strings to their reverses, up to absolute value.
Both {\bf s} and {\bf t} are nonempty and contain an even number of tiles. Furthermore,  {\bf s} and {\bf t} both have odd length and so are not equal to all of {\bf c}.  If ${\bf s}={\bf t}$, then  there exists a pure seam either before {\bf s}, after {\bf s}, or both, a contradiction.  Thus ${\bf s}\ne {\bf t}$. 

Assume first that {\bf t} is strictly shorter than {\bf s}. Then {\bf t} lies in the interior of {\bf s} and therefore in the interior of {\bf c}.
Let $D_1$ and $D_2$ be the {\bf b}-connectors that immediately precede and follow {\bf t}, respectively.  If $D_1=(0)$, then $D_1$ lies in {\bf s} and $\sigma_{\bf s}(D_1)=(0)$ is in {\bf s}. Thus $D_2$ begins with zero and so must also be $(0)$, a contradiction, since we can now extend {\bf t} by at least one more {\bf b}-connector and {\bf b}-tile at each end.  Similarly, $D_2$ is not zero. Thus both $D_1$ and $D_2$ are maximal components in {\bf c} and hence, if either lies in {\bf s}, is a maximal component in {\bf s}. The involution $\sigma_{\bf s}$ must take maximal components in {\bf s} to maximal components in {\bf s}. So  if $D_1$ and $D_2$ both lie  in {\bf s}, then  $\sigma_{\bf s}(D_1)=\pm D_2$, again contradicting the fact that {\bf t} is maximal.  Next, suppose that neither $D_1$ nor $D_2$ lie in {\bf s}. 
In this case it follows that {\bf s} must also be preceded and followed by {\bf a}-connectors, which we call $C_1$ and $C_2$, respectively. But now $C_i \cap D_i$ is nonempty for $i=1,2$, and so we  have that each $C_i=(0)$ and is centered in $D_i$. But now we see that {\bf s} is not maximal. The final case is that one of $D_1$ or $D_2$  lies in {\bf s} while the other does not. Suppose $D_1$ lies in {\bf s}. If it lies in the interior of {\bf s} then $\sigma_{\bf s}(D_1)$ lies in the interior of {\bf s} and so is a  maximal component in {\bf s}. Thus we have that $\sigma_{\bf s}(D_1)=\pm D_2$, a contradiction. 
So we may assume that ${\bf s}=(D_1, {\bf t}, \pm \sigma_{\bf s}(D_1))$. Since $D_2$ does not lie in {\bf s}, an {\bf a}-connector $C_2$ follows {\bf s} and so must intersect $D_2$. Thus $C_2=(0)$ and is centered in $D_2=\pm (\sigma_{\bf s}(D_1), 0, \sigma_{\bf s}(D_1))$.  
 Let {\bf x} and {\bf y} be those parts of {\bf c} that precede {\bf s} and follow $C_2$, respectively. If we concatenate {\bf x} and {\bf y} we obtain a shorter element of ${\cal S}_{even}$ that parses with respect to {\bf a} and {\bf b}, a contradiction.  
Similar arguments work if $D_1$ is not contained in {\bf s}  and $D_2$ is contained in {\bf s} 
or if {\bf s} is shorter than {\bf t}.  We have now proven the second property.

Since {\bf a} and {\bf b}-connectors do not intersect, every {\bf a}-connector lies in a {\bf b}-tile.  Suppose an {\bf a}-connector $C$ is not a maximal component  in the {\bf b}-tile. Then $C$ is a proper subset of a maximal component $D$  in the {\bf b}-tile. If $C$ is not zero, then a zero immediately precedes or follows $C$. But this is impossible because  {\bf a}-tiles do not begin or end with zero. If $C$ is zero,  then $D=(D_1, 0, D_2)$. Suppose $C$ is preceded in the {\bf a}-parsing by the {\bf a}-tile $A$ and followed by the {\bf a}-tile $\pm A^{-1}$.  Now $A$ ends with $D_1$ which is  a maximal component in $A$ and $A^{-1}$ begins with $D_2$ which is  a maximal component in $A^{-1}$. Hence $D_1=\pm D_2$ and $C=(0)$ is centered in $D$. This completes the proof of the third property. A similar argument gives the last property. \hfill $\square$

\begin{figure}[t]
    \begin{center}
    \leavevmode
    \scalebox{.30}{\includegraphics{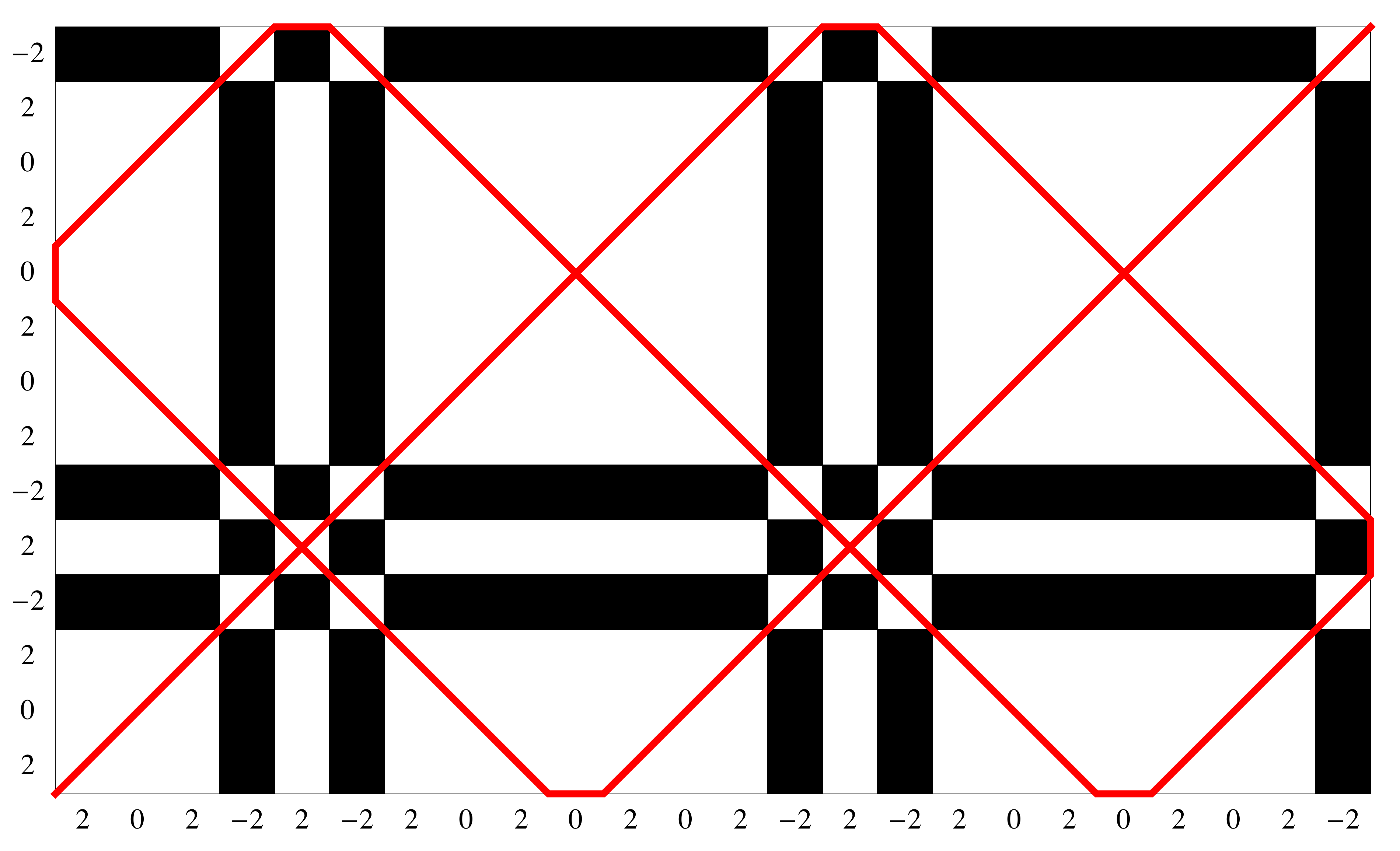}}
    \end{center}
\caption{The product ${\bf a} \times {\bf b}$ colored with respect to maximal component decompositions. All signs are $(+,+)$ and are omitted.} 
\label{cPath4}
\end{figure}

Before continuing, we remark further on the geometric interpretation of {\bf c} as a path in ${\bf a} \times {\bf b}$ as introduced in Figure~\ref{cPath1}.  Given {\bf a} and {\bf b}, write  ${\bf a}=(C_1, C_2, \dots, C_{r-1})$ and  ${\bf b}=(D_1, D_2, \dots, D_{s-1})$, where each $C_i$ and $D_j$ is a maximal component in {\bf a} and {\bf b}, respectively. 
We now form the product ${\bf a}\times {\bf b}$ and color it checkerboard fashion with respect to the maximal component decompositions by shading the cell $C_i \times D_j$ white  if $i+j$ is even, and black otherwise.  In Figure~\ref{cPath4} we give an example of {\bf a}, {\bf b}, and {\bf c} such that {\bf c} parses with respect to both {\bf a} and {\bf b}. As already described, each diagonal segment of {\bf c} can be labeled with a pair of signs $(\epsilon, \eta)$ indicating that this segment of {\bf c} lies in an {\bf a}-tile of sign $\epsilon$ and a {\bf b}-tile of sign $\eta$. Notice that in this example, {\bf c} can only pass through an {\it interior} white region $C_i \times D_j$ with $1<i<r-1$ and $1<j<s-1$ along its diagonal and hence if and only if $\epsilon C_i= \eta D_j$ where that segment of {\bf c} is labeled $(\epsilon, \eta)$. Similarly, {\bf c}  only passes through an {\it exterior} white region, such as $C_1 \times D_j$, if either $ \epsilon C_1= \eta D_j$ or the length of $D_j$ is one more than twice the length of $C_1$ (and the signs are appropriate). We can express this last condition more succinctly by thinking of a connector as both a vector and an even number equal to the sum of its entries.  The condition now becomes $2 \epsilon C_1= \eta D_j$. For example, in Figure~\ref{cPath4}, {\bf c} passes through $C_1 \times D_5$ and $2C_1=D_5=8$. However it may also pass through a square exterior region like  $C_6 \times D_2$ with  $C_6=D_2=-2$. When the path reaches the boundary of ${\bf a} \times {\bf b}$ it traverses a connector, which by Lemma~\ref{properties of c}, is either the entire boundary of a black region, or the central zero of a white region. 

\begin{lemma}\label{path exists}Let ${\bf a}, {\bf b}, {\bf c} \in S_{even}$ and suppose that {\bf c} parses with respect to both ${\bf a}$ and ${\bf b}$ and that no shorter vector in $S_{even}$ has this property. Then {\bf c} can be represented by a path in ${\bf a} \times {\bf b}$ of the type described in Figure~\ref{cPath4}. 
\end{lemma}
{\bf Proof:} Form the product ${\bf a}\times {\bf b}$ and initially think of the rows and columns as labeled with the individual entries of {\bf a} and {\bf b}, rather than maximal components. We now define a map $f:{\bf c} \to {\bf a}\times {\bf b}$ taking each entry of {\bf c} to an oriented line segment as follows.
If $x$ is an entry of {\bf c}, then by Lemma~\ref{properties of c}, $x$ is either in both an {\bf a}-tile and a {\bf b}-tile, or in an {\bf a}-tile and a {\bf b}-connector, or in an {\bf a}-connector and a {\bf b}-tile. In the first case, $x$ determines a specific row and column in the product in an obvious way. If $x$ is in the tiles $\epsilon {\bf a}^\alpha$ and $\eta {\bf b}^\beta$, where $\epsilon, \eta, \alpha$ and $\beta$ are in $\{-1, +1\}$, then $f(x)$ is defined to be the oriented line segment in that row and column whose direction is given by the vector $\left <\alpha, \beta\right >$. Moreover, this diagonal segment if labeled $(\epsilon, \eta)$. 
In the second  case, $x$ lies in a specific row determined by its place in the {\bf a}-tile, but since it does not lie in the {\bf b}-tile, it does not correspond to any column. Instead, it corresponds to  a place on the left or right boundary of the product. In particular, if $x$ lies in a {\bf b}-connector that follows the tile $\eta {\bf b}$, then $x$ lies on the right boundary of the product, and if it lies in a {\bf b}-connector that follows the tile  $\eta {\bf b}^{-1}$, then $x$ lies on the left boundary. Finally, if $x$ lies in $\epsilon {\bf a}^\alpha$, then the orientation of $f(x)$ is given by the vector $\left < 0, \alpha\right >$. Similarly, the third case gives oriented segments on the top or bottom edge of ${\bf a}\times {\bf b}$. It is now a simple matter to verify that the image of {\bf c} under $f$ is a connected, coherently oriented path which never travels through the interior of a black region in the checkerboard coloring. Further, consecutive diagonal segments are labeled with the same pair of signs. It is also important to note that each segment of the path (except the first segment) uniquely determines the preceding segment. Hence the path cannot traverse any given segment more than once.  \hfill $\square$

We now continue with the assumption that {\bf a}, {\bf b} and {\bf c} satisfy the hypothesis of Lemma~\ref{path exists}. Thus we may represent {\bf c} as a path in ${\bf a}\times{\bf b}$ as  described in Figure~\ref{cPath4}. The following lemma draws conclusions about the {\bf a} and {\bf b}-connectors that lie along the edges of the product ${\bf a}\times{\bf b}$.  

\begin{lemma}\label{properties of connectors}Suppose ${\bf a}, {\bf b}, {\bf c} \in {\cal S}_{even}$,  {\bf c} parses with respect to both ${\bf a}$ and ${\bf b}$, and that no shorter vector in ${\cal S}_{even}$ has this property. Then the path {\bf c} in ${\bf a}\times {\bf b}$ satisfies all of the following.
\vspace{-.25 in}
\begin{enumerate}
\item The connectors on the union of the bottom  and  left edges of the product ${\bf a}\times {\bf b}$ are either all zero or all nonzero.
\item The connectors on the union of the  top and  right edges of the product ${\bf a}\times {\bf b}$ are either all zero or all nonzero.
\item No two {\bf a}-connectors lie in the same maximal component of a {\bf b}-tile.
\item No two {\bf b}-connectors lie in the same maximal component of an {\bf a}-tile.
\item Every diagonal segment of {\bf c} is labeled $(+,+)$ or $(-,-)$.
\item Either {\bf c} has no mixed seams, or exactly two mixed seams. 
\end{enumerate}
\end{lemma}

\noindent{\bf Proof:} The proof of this lemma is quite long and technical. We begin with the following claim.

\noindent{\bf Claim:} All connectors along a single edge of the product are either all zero or all nonzero.
 
Suppose instead that a pair of connectors exist on one edge with one being zero and the other non zero. We may assume the two connectors are adjacent, that is, no other connectors on this edge lie between them. For concreteness, suppose they are two {\bf a}-connectors that lie on the the boundaries of the cells $C_{r-1} \times D_i$ and $C_{r-1} \times D_j$, respectively, with $i<j$, and that the {\bf a}-connector in $C_{r-1} \times D_i$ is zero. The argument we give will apply equally well to the other possible cases. We know that $C_{r-1} \times D_i$ must be a white cell and that $C_{r-1} \times D_j$ must be black, thus $i$ is even and $j$ is odd. 
We focus now on the gap between them,  $k=j-i-1$,  which is even.
Since $D_i$ has a central zero, we may write $|D_i|=2|C_{r-1}|$. Finally, let $\alpha$ be the diagonal segment of {\bf c} with slope $-1$ which leaves the cell $C_{r-1} \times D_i$ headed down and to the right, and $\beta$ the diagonal segment of {\bf c} with slope $+1$ which leaves the cell $C_{r-1} \times D_j$ headed down and to the left.

\noindent{\bf Case A}, $k\ge 2r$:
In this case the segment $\alpha$ will bounce off the bottom edge of the product and return to the top edge giving an {\bf a}-connector between the adjacent {\bf a}-connectors in positions $i$ and $j$, a contradiction. In particular, $\alpha$ will either arrive at a zero {\bf a}-connector on the bottom edge at $C_1 \times D_{i+r-2}$ or at  a nonzero {\bf a}-connector  at $C_1 \times D_{i+r-1}$. The path will then leave the {\bf a}-connector on the bottom edge and reach an {\bf a}-connector on the top edge at $C_{r-1}\times D_{i+2r-x}$, where $x\in\{1,2,3,4\}$. In any of these cases, $i+2r-x<j$ since $k\ge 2r$.

\noindent{\bf Case B}, $k=2r-2$:
With $k=2r-2$ we have that $i+2r-1\le j$. Thus the argument of the previous case applies unless $\alpha$  and $\beta$ are joined together at a nonzero {\bf a}-connector on the bottom edge at $C_1\times D_{i+r-1}$. In this case, consider the column of the product given by ${\bf a}\times D_{i+r-1}$. The path {\bf c} must cross this column an odd number of times greater than one, and hence some other part of {\bf c} which crosses through this column must give an {\bf a}-connector on the top edge of the product between the adjacent {\bf a}-connectors in positions $i$ and $j$, a contradiction.

\noindent{\bf Case C}, $k<r-1$:
With the gap small, we cannot argue that another {\bf a}-connector lies on the top edge between the adjacent {\bf a}-connectors at positions $i$ and $j$. Instead we will reach a contradiction in another way. Because $k<r-1$, $\alpha$ passes through $C_{r-1-k} \times D_{j-1}$ which is in the second row, or higher. Similarly, the path $\beta$ passes through $C_{r-1-k} \times D_{i}$. Remember now that if {\bf c} passes diagonally through the cell $C_u \times D_v$ then $|C_u|=|D_v|$. Thus  we have 
$$|C_{r-1}|=|D_{j-1}|=|C_{r-1-k}|=|D_i|=2| C_{r-1}|,$$ a contradiction, since maximal components are never zero.

\noindent{\bf Case D}, $r-1<k<2r-2$:
In this case, $\alpha$ and $\beta$ will cross each other before reaching the bottom edge unless $k=2r-4$ and $\alpha$ and $\beta$ are joined at a zero {\bf a}-connector at $C_1 \times D_{i+r-2}$. If this happens, then as in Case B, some other part of {\bf c} crosses the column  ${\bf a}\times D_{i+r-2}$ and produces an {\bf a}-connector on the top edge of the product between the adjacent {\bf a}-connectors in positions $i$ and $j$, a contradiction. If $\alpha$ and $\beta$ are not joined at a zero {\bf a}-connector on the bottom edge, then they cross and reach two separate {\bf a}-connectors on the bottom edge separated by a gap less than $r-2$. Hence by Case C, these two {\bf a}-connectors on the bottom edge are either both zero or both nonzero. Consider first the case where they are both zero. We now have that $\alpha$ will reach the zero {\bf a}-connector at $C_1 \times D_{i+r-2}$ and then travel back up passing through $C_{k-r+3}\times D_{j-1}$. Moreover, $\beta$ will reach the zero {\bf a}-connector at $C_1 \times D_{j-r+1}$ and then travel back up passing through $C_{k-r+3}\times D_{i}$. Thus
 $$|C_{r-1}|=|D_{j-1}|=|C_{k-r+3}|=|D_i|=2 |C_{r-1}|,$$ 
 a contradiction. If instead, the two connectors on the bottom edge are both nonzero, this lowers each of the rebounding paths by two units so they now travel through  $C_{k-r+1} \times D_{i}$ and  $C_{k-r+1} \times D_{j-1}$. We now have  the contradiction
  $$|C_{r-1}|=|D_{j-1}|=|C_{k-r+1}|=|D_i|=2 |C_{r-1}|.$$ 
  
\noindent{\bf Case E}, $k=r-1$:
We now have that $\alpha$ and $\beta$ cross and reach two separate {\bf a}-connectors on the bottom edge separated by a gap less than $r-1$. Hence by Case C, these two {\bf a}-connectors on the bottom edge are either both zero or both nonzero. If they are both zero, we are in the situation of Case D. If they are both nonzero, then they lie at 
$C_1 \times D_{i}$ and $C_1 \times D_{j-1}$. 
Now $\alpha$ and $\beta$ cross each other in $C_{\frac{r-1}{2}} \times D_{i+\frac{r-1}{2}}.$ The path {\bf c} must cross column ${\bf a}\times D_{i+\frac{r-1}{2}}$ an odd number of times. If it crosses this column above row $\frac{r-1}{2}$ it will extend up to the top edge and create another {\bf a}-connector between columns $i$ and $j$, a contradiction. Thus it must cross the column  below row $\frac{r-1}{2}$, and create another {\bf a}-connector on the bottom edge of the product between the two nonzero {\bf a}-connectors in $C_1 \times D_i$ and $C_1 \times D_{j-1}$. This {\bf a}-connector must be nonzero because it is less than $r-1$ away from a nonzero {\bf a}-connector. Suppose it is located in $C_1 \times D_l$. The diagonal part of {\bf c} leaving the {\bf a}-connector at $C_1 \times D_l$ headed up and to the left passes through $C_{l-i} \times D_i$.  
The diagonal part of {\bf c} leaving the {\bf a}-connector at $C_1 \times D_l$ headed up and to the right passes through $C_{j-l-1} \times D_{j-1}$. Furthermore, $\alpha$ and $\beta$ pass though column ${\bf a}\times D_l$ at $C_{j-l-1} \times D_{l}$ and $C_{l-i} \times D_{l}$, respectively. Thus we have 
 $$|C_{r-1}|=|D_{j-1}|=|C_{j-l-1}|=|D_l|=|C_{l-i}|=|D_i|=2 |C_{r-1}|,$$ a contradiction.

 \noindent{\bf Proof of 1 and 2:} Now that we know all connectors along a single edge of the product are either all zero or all nonzero, we address the first and second properties of the Lemma. Consider the first {\bf a}-connector on the bottom edge (that is, closest to the lower left corner) and the first {\bf b}-connector on the left edge (again, closest to the lower left corner.) These must be connected by a diagonal part of {\bf c}. We will show that they are either both zero or both nonzero. A similar argument applies near the upper right corner of the product.

Suppose a zero {\bf b}-connector is located at $C_i \times D_1$ and is connected to a nonzero {\bf a}-connector at $C_1 \times D_{i+1}$. Note that $D_i=C_i$ for $0<i<r-1$ since the initial diagonal part of {\bf c} passes through $C_i \times D_i$ and is labeled $(+,+)$. Thus $D_1=C_1=|D_i|=|C_i|=2|D_1|$, a contradiction. A similar argument shows that   the {\bf a}-connector cannot be zero and the {\bf b}-connector nonzero. This completes the proof of the first two properties.

 \noindent{\bf Proof of 3 and 4:} Suppose now there are a pair of {\bf a}-connectors that lie in the same maximal component $D_j$ of the {\bf b}-tile. (The case of two {\bf b}-connectors lying in the same maximal component of the {\bf a}-tile is similar.) By Lemma~\ref{properties of c}, each {\bf a}-connector is either all of $D_j$, or the central zero of $D_j$. Because $r$ is odd, the two cells $C_1 \times D_j$ and $C_{r-1}\times D_j$ have opposite colors. Thus if the two {\bf a}-connectors both lie on the same edge of the product, they would have to either both be the central zero of $D_j$, or both be all of $D_j$. But these cases are impossible, as they would require the path  to traverse the same segment twice. We ruled this out in the proof of Lemma~\ref{path exists}. Thus one {\bf a}-connector lies on the top edge, and one on the bottom edge.

Without loss of generality, assume that there is a zero {\bf a}-connector on the top edge of the product and a nonzero {\bf a}-connector on the bottom edge. But now all connectors on the top and right edges are zero and all the connectors on the left and bottom edges are non zero. If we now remove the last row and last column from the product, we obtain a new product ${\bf a}'\times {\bf b}'$ of the truncated tiles ${\bf a}'$ and ${\bf b}'$ obtained from {\bf a} and {\bf b} by removing the last maximal components  $C_{r-1}$ and $D_{s-1}$, respectively. Notice that in this new product we may obtain a path ${\bf c}'$ from {\bf c} by replacing what had been zero connectors along the top and right edges with nonzero connectors in the obvious way. The new path has no nonzero connectors, but still has at least one pair of ${\bf a}'$-connectors that lie in the same maximal component of the ${\bf b}'$-tile. 
Since none of the connectors, in either parsing of ${\bf c}'$ are zero, it is not hard to show by induction that all of the ${\bf a}'$-connectors along the top edge of the product occur at positions $C_{r-2}\times D_j$ where $j$ is an odd multiple of $\gcd(r-1, s-1)$ and along the bottom edge at positions $C_{1}\times D_j$ where $j$ is an even multiple of $\gcd(r-1, s-1)$. Therefore, two ${\bf a}'$-connectors cannot lie in the same maximal component of the ${\bf b}'$-tile.

\noindent{\bf Proof of 5:} We now wish to show that each diagonal segment of the path {\bf c} is labeled either $(+,+)$ or $(-,-)$. As already mentioned, if two diagonal segments of {\bf c} intersect, then the product of the signs labeling each segment is the same. Thus if the set of diagonal segments form a connected subset of the plane, as they do in Figure~\ref{cPath1}, it follows that all diagonal segments are labeled $(+,+)$ or $(-,-)$ since the first segment is labeled $(+,+)$. But as Figure~\ref{cPath4} illustrates, the set of diagonal segments may not form a connected set. If the path {\bf c} does not pass through the upper left corner, then the left-most diagonal segment of slope $+1$ will not intersect any other diagonal segment. Similarly, the right-most diagonal segment of slope $+1$ might also be disconnected  from the other diagonal segments. We need to show that in this case, the isolated segments are still labeled $(+,+)$ or $(-,-)$. 

If two consecutive diagonal segments are separated by an {\bf a}-connector, then they lie in the same {\bf b}-tile and hence the second sign labeling each segment must be  the same. But they lie in separate {\bf a}-tiles and so the first signs labeling each segment might be different. However, if the {\bf a}-connector is zero, we know that it is the central zero of a maximal component  in the {\bf b}-tile and this forces the signs of the {\bf a}-tiles on either side of this {\bf a}-connector to be the same. Thus the pair of signs labeling a diagonal segment can only change at a nonzero connector.  It suffices to consider the case of an isolated diagonal segment of slope $+1$ with nonzero connectors at both ends.  This implies  that all {\bf a} and {\bf b}-connectors are nonzero. As we have already seen,  
 the {\bf a}-connectors along the top edge are located at odd multiples of $d$ and along the bottom edge are located at even multiples of  $d$,  where $d=\gcd (r, s)$.  Similarly, along the left edge, {\bf b}-connectors are located at even multiples of $d$, and along the right edge at odd multiples of $d$. Hence the left-most diagonal segment of {\bf c} with slope $+1$, which we will call $\Sigma$,  connects the {\bf b}-connector at $C_{(\frac{r}{d}-1)d}\times D_1$ to the {\bf a}-connector at $C_{r-1}\times D_d$. The next  segment with slope $+1$ connects the {\bf b}-connector at $C_{(\frac{r}{d}-3)d}\times D_1$ to the {\bf a}-connector at $C_{r-1}\times D_{3d}$. Intersecting this segment are two diagonal segments of slope $-1$: one has the {\bf a}-connector at $C_{r-1} \times D_d$ at its end, the other has the {\bf b}-connector at $C_{(\frac{r}{d}-1)d} \times D_1$ at its end. These four diagonal segments form a square. Note that except possibly for $\Sigma$, each of these diagonal segments is labeled either $(+,+)$ or $(-,-)$. Suppose $\Sigma$ is labeled $(\epsilon, \eta)$. Finally we have that
$$\epsilon \eta C_{r-1}= D_{d-1}=C_{r-2d+1}=D_{d+1}=C_{r-1}.$$
Hence $\epsilon \eta=1$ as desired. 

\noindent{\bf Proof of 6:} If {\bf c} has mixed seams, then this corresponds to the path {\bf c} passing through the upper left, or lower right corner of the product ${\bf a}\times {\bf b}$. If, for example, the path passes through the upper left corner, then there are nonzero connectors on the top and left edges of the product, and hence on all edges of the product. Similarly, if {\bf c} passes through the lower right corner, we again find that all connectors are nonzero. As mentioned already, we now know the locations of the connectors in terms of $d=\gcd(r,s)$. Having a mixed seam at either corner now implies that $d=1$, which in turn implies there are mixed seams at both corners. \hfill $\square$

With these lemmas established we now proceed to prove our main theorem.

{\bf Proof of Theorem~\ref{main theorem}:} Suppose that $K_1=\Phi(\hat{\bf a})$ and $K_2=\Phi(\hat{\bf b})$ are incomparable knots with an upper bound $K=\Phi({\hat{\bf c}})$. Then there are representatives of $\hat{\bf a}, \hat{\bf b}$, and $\hat{\bf c}$, which we call {\bf a}, {\bf b}, and {\bf c}, respectively, such that {\bf c} parses with respect to both {\bf a} and {\bf b}. Furthermore, we may assume that {\bf c} is shortest among all such {\bf c} and that {\bf a} is shorter than {\bf b}. Decompose both {\bf a} and {\bf b} into  maximal components and consider {\bf c} as a path in ${\bf a} \times {\bf b}$. Suppose an {\bf a}-connector is zero and is the central zero of the maximal component $D_j$ in the {\bf b}-tile. Split $D_j$ into three components $(D_j', 0, D_j')$. Similarly, split any maximal component in the {\bf a}-tile whose central zero is a {\bf b}-connector into three components. This gives a new decomposition of {\bf a} and {\bf b} into components, not necessarily maximal ones. We leave it to the reader to show that the conclusions of Lemma~\ref{properties of connectors} now imply that these two decompositions induce a decomposition of the path {\bf c} into a sequence of components that satisfy the hypothesis of Proposition~\ref{gcd result}. We now obtain the conclusion of Theorem~\ref{main theorem}. \hfill $\square$

Using Theorem~\ref{main theorem}, we can easily decide if $\{K_1, K_2\}$ has an upper bound. Moreover, given a single 2-bridge knot $K_1$, we can now classify all other 2-bridge knots $K_2$ such that $\{K_1, K_2\}$ has an upper bound. In the following example we do this for the  figure-eight  and trefoil knots.

{\bf Example 3:} Consider the figure eight knot $K_1$ represented by ${\bf a}=(2,2)$. Suppose $K_1$ and $K_2=\Phi(\hat{{\bf b}})$ have an upper bound, but $K_1$ and $K_2$ are not comparable. Then by Theorem~\ref{main theorem}, $(2,2)=({\bf w}^p,{\bf e})$ which implies that {\bf e} is empty, $n=m=2$ and $p=1$.  Hence  ${\bf b}=(2, 2, 2, \dots, 2)$, a vector of $2 n$ twos with $n \ne 1 \ (\mbox{mod } 3)$. The last condition is necessary to ensure that $K_2$ is not greater than $K_1$. Similarly, the knots that share an upper bound with the trefoil, but are not comparable to the trefoil,  are represented by  ${\bf b}=(2, -2, 2, \dots, -2)$, a vector of $2 n$ twos, alternating in sign, with $n \ne 1 \ (\mbox{mod } 3)$.

If  $K_1$ and $K_2$ are represented by the vectors {\bf a} and {\bf b}, respectively, then we can produce all possible shortest-length vectors {\bf c} that parse with respect to both {\bf a} and {\bf b}. Theorem~\ref{main theorem} gives us one such vector {\bf c}, and if {\bf c} has no mixed seams, it is unique. However, if {\bf c} contains two mixed seams, then by negating the entries (as a group) that appear  between, or after the seams, we obtain all other shortest-length vectors that parse with respect to both {\bf a} and {\bf b}. It is not hard to show that a shortest-length vector that parses with respect to both {\bf a} and {\bf b} represents a least upper bound of $\{K_1, K_2\}$. However, not every least upper bound of $\{K_1, K_2\}$ is represented by a shortest-length vector. For example, let {\bf a}, {\bf b}, and {\bf c} be as in Example~1. Let ${\bf c}'$ be obtained from {\bf c} by negating all the entries between the two mixed seams. Finally, let ${\bf d}=(c, -2,0,-2, c'^{-1}, 2, c')$. Clearly {\bf d} parses with respect to both {\bf a} and {\bf b}, but we leave it to the reader to check that {\bf d} does not parse with respect to any other vector. Thus {\bf d} represents a least upper bound of $\{K_1, K_2\}$ and yet is not of minimal length.

We now turn our attention to determining if any given set of 2-bridge knots has an upper bound. We need the following theorem.

\begin{theorem} \label{std form}Suppose ${\bf a} \in {\cal S}_{even}$ can be expressed as 
\begin{equation}\label{special form}{\bf a}=(({\bf e}, m, {\bf e}^{-1}, n)^q, {\bf e}),\end{equation}
for some ${\bf e}\in {\cal S}_{even}$, possibly empty, $m$ and $n$ even integers, and $q$ a natural number. Then $m$ and $n$ are unique. Moreover, if ${\bf e}_0$ is the shortest possible choice for {\bf e}, with $q_0$ the associated power, then  ${\bf e}_0$ is unique and every other possible choice of {\bf e} is given by ${\bf e}=(({\bf e}_0, m, {\bf e}_0^{-1}, n)^l, {\bf e}_0)$ with $(2l+1)q+l=q_0$.
\end{theorem}
\noindent{\bf Proof:} If {\bf a} can be expresses as in Equation~\ref{special form}, then clearly there is a unique vector ${\bf e}_0$ of shortest length for which this is true. Suppose first that ${\bf e}_0$ is empty so that ${\bf a}=(m,n)^{q_0}$. Since {\bf a} does not begin or end with zero, neither $m$ nor $n$ are zero. Hence this form represents the unique  decomposition of {\bf a} into maximal components. Therefore, if {\bf a} can also be written as ${\bf a}=(i,j)^r$, then $i=m, j=n$ and $r=q_0$. Suppose now that {\bf a} can also be expressed as ${\bf a}=(({\bf e}, i, {\bf e}^{-1}, j)^q, {\bf e})$ and ${\bf e}=(C_1, C_2, \dots, C_{2l})$ is the unique decomposition of {\bf e} into maximal components. Since the last maximal component of {\bf a} is both $n$ and $C_{2l}$, we see that $C_{2l}=n$. If $i=0$ then the $2l$-th maximal component of {\bf a} is both $n$ and $(C_{2l}, 0, C_{2l})=2n$, a contradiction. Similarly, we cannot have $j=0$. But now both  $i$ and $j$ are maximal components, and the uniqueness of maximal component decompositions gives that $i=m, j=n$ and ${\bf e}=(m,n)^l$ with $(2l+1)q+l=q_0$.

Next consider  ${\bf a}=(({\bf e}_0, m, {\bf e}_0^{-1}, n)^{q_0}, {\bf e}_0)$ where ${\bf e}_0$ is nonempty. Furthermore, assume that ${\bf e}_0$ is the shortest possible such vector, and also that ${\bf a}=(({\bf e}, i, {\bf e}^{-1}, j)^q, {\bf e})$ for some other vector {\bf e} in ${\cal S}_{even}$ which is longer than ${\bf e}_0$. We now have that {\bf a} parses with respect to both ${\bf e}_0$ and {\bf e}. We first consider the case where {\bf e} parses with respect to ${\bf e}_0$. We now have two parsings of {\bf a} with respect to ${\bf e}_0$: the one given by ${\bf a}=(({\bf e}_0, m, {\bf e}_0^{-1}, n)^{q_0}, {\bf e}_0)$ and the one given by first parsing {\bf a} with respect to {\bf e} and then parsing {\bf e} with respect to ${\bf e}_0$. Since these must be the same, we have that $i=m, j=n$ and ${\bf e}=(({\bf e}_0, m, {\bf e}_0^{-1}, n)^l, {\bf e}_0)$ with $(2l+1)q+l=q_0$. If {\bf e} does not parse with respect to ${\bf e}_0$, then by Theorem~\ref{main theorem}, we must have ${\bf e}_0=(({\bf g}, r, {\bf g}^{-1}, s)^t, {\bf g})$ and ${\bf e}=(({\bf g}, u, {\bf g}^{-1}, v)^w, {\bf g})$  for some vector {\bf g}, possibly empty, $r, s, u$, and $v$ even, and $t$ and $w$ natural numbers. Suppose  first that  {\bf g} is empty. It now follows that none of $r, s, u$, or $v$ are zero. Comparing maximal component decompositions of {\bf a} now gives that $u=r$ and $v=s$. Furthermore,  using an argument similar to one already given, if $m=0$ we obtain the contradiction $s=2s$ and if $n=0$ the contradiction $r=2r$. We now obtain $r=m$ and $s=n$, and contradict the assumption that {\bf a} cannot be expressed as in Equation~\ref{special form} with ${\bf e}_0$ empty. Hence {\bf g} is non empty. In this case, we can think of {\bf a} as parsing with respect to {\bf g} in two ways: first with respect to ${\bf e}_0$ and then with respect to {\bf g}, or first with respect to {\bf e} and then with respect to {\bf g}. These two points of view imply that $u=r=m$ and $v=s=n$. But now {\bf g} is shorter than ${\bf e}_0$, a contradiction.
We conclude that the case of  {\bf e} not parsing with respect to ${\bf e}_0$ is impossible. \hfill $\square$

We are now in a position to determine if any set $S$ of 2-bridge knots has an upper bound. As already mentioned, this is not possible if $S$ is infinite. So assume $S=\{K_1, K_2, \dots, K_l\}$. Now if any one of these knots is smaller than another, we can remove it from the set and it suffices to find an upper bound for the remaining set of knots. Hence we may assume the knots in $S$ are pairwise incomparable. If $S$ has an upper bound, then each pair of knots in $S$ has an upper bound. In particular, comparing $K_1$ to each of the other knots $K_i$ implies that $K_1$ and $K_i$ are represented by the vectors $(({\bf e}_i, m_i, {\bf e}_i^{-1}, n_i)^{p_i}, {\bf e}_i)$ and $(({\bf e}_i, m_i, {\bf e}_i^{-1}, n_i)^{q_i}, {\bf e}_i)$, respectively, for each $i$. But now we have represented $K_1$ in $l-1$ different ways. By Theorem~\ref{std form}, it follows that $m_i=m_2$ and $n_i=n_2$ for all $i$. Let $m=m_2$ and let $n=n_2$. Moreover, there is a unique vector  ${\bf e}_0$ and exponent $p_0$ such that ${\bf e}_i=(({\bf e}_0, m, {\bf e}_0^{-1}, n)^{r_i}, {\bf e}_0)$ where $(2 r_i+1)p_i+r_i=p_0$ for all $i$. But now $K_1$ is represented by the vector $(({\bf e}_0, m, {\bf e}_0^{-1}, n)^{p_0}, {\bf e}_0)$ and each knot $K_i$ with $i>1$ is represented by the vector $(({\bf e}_0, m, {\bf e}_0^{-1}, n)^{(2 r_i+1)q_i+r_i}, {\bf e}_0)$. Hence an upper bound for the entire set of knots exists and is represented by the vector $(({\bf e}_0, m, {\bf e}_0^{-1}, n)^Q, {\bf e}_0)$, where $2Q+1$ is the least common multiple of $\{2p_0+1, (2r_2+1)(2q_2+1),   (2r_3+1)(2q_3+1), \dots,  (2r_l+1)(2q_l+1)\}$. Notice that in fact what we have proven is the following generalization of Theorem~\ref{main theorem}.

\begin{theorem}\label{second main theorem}The 2-bridge knots $K_1=\Phi(\hat{\bf a}_1), K_2=\Phi(\hat{\bf a}_2),\dots ,K_l=\Phi(\hat{\bf a}_l)$ are pairwise incomparable and $\{K_1, K_2, \dots, K_l\}$ has an upper bound with respect to the Ohtsuki-Riley-Sakuma partial order if and only if for each $i$
$${\bf a}_i=({\bf w}^{p_i}, {\bf e}),$$
where $\bf e$ is some (possibly empty) vector in ${\cal S}_{even}$,  ${\bf w}=({\bf e}, {\bf m}, {\bf e}^{-1}, {\bf n})$, $m$ and $n$ are even integers, and finally, $2p_i+1$ does not divide  $2p_j+1$ for all $i\ne j$.
\end{theorem}

Theorem~\ref{second main theorem} has the following interesting corollary.

\begin{corollary} The set of 2-bridge knots $\{K_1, K_2, \dots, K_l\}$ has an upper bound with respect to the Ohtsuki-Riley-Sakuma partial order if and only if each pair $\{K_i, K_j\}$ has an upper bound.
\end{corollary}

\bibliographystyle{plain}
\vskip -.275 in

 \end{document}